\documentclass[12pt]{article}
\usepackage{amssymb,amsmath}
\usepackage[dvips]{graphicx}
\def\dH{{\mathbb{H}}}
\def\dR{{\mathbb{R}}}

\def\dC{{\mathbb{C}}}

\def\dM{{\mathbb{M}}}
\def\ld{{\lambda}}

\def\szs{{\mathcal{S}}}

\def\conv{{\mbox{conv}}}
\def\span{{\mbox{span}}}

\def\diag{{\mbox{diag}}}

\def\szs{{\\[.5cm]}}
\newcommand{\beq}[2]{\begin{equation}\label{#1}#2}

\begin{document}
\noindent
\textbf{Normal matrix compressions}
\szs
6 December 2011
\szs
John Holbrook, Nishan Mudalige, Rajesh Pereira
\szs
\textbf{Abstract:} The recently developed theory of higher--rank numerical ranges originated in
problems of error correction in quantum information theory but its mathematical implications now include
a quite satisfactory understanding of \emph{scalar} compressions of complex matrices. Here our aim is to
make some first steps in the more general program of understanding \emph{normal} compressions. We establish
some general principles for the program and make a detailed study of rank--two normal compressions.
\szs
\szs
\szs
\textbf{1: Introduction}
\szs
Given a linear operator $T$ on a complex Hilbert space $\dH$, and any orthogonal projection $P$, we say that
$PT|_{P\dH}$ is a \textbf{compression} of $T$. If $\dH=\dC^N$ and $T$ is represented by a matrix $M\in\dM_N$
(the $N\times N$ complex matrices), a second matrix $C$ represents
a compression of $T$ (or a compression of $M$) iff there is a unitary matrix $U$ such that $C$ is a NW corner of $UMU^*$.
If $C$ is $k\times k$ we say it is a rank--$k$ compression of $M$. There is a rich history of results that allow us to
identify compressions by means of intrinsic criteria. A classic example is the Cauchy interlacing theorem [Cau], along with its
converse [FP], which may be expressed as follows.
\szs
\textbf{Theorem 1:} If $M\in\dM_N$ is Hermitian, with eigenvalues
\[
a_1\leq a_2\leq\dots\leq a_N,
\]
then $C$ is a rank--$k$ compression of $M$ iff $C$ is Hermitian with eigenvalues $b_j$ satisfying
\[
a_1\leq b_1\leq a_{N-k+1}, a_2\leq b_2\leq a_{N-k+2},\dots,a_k\leq b_k\leq a_N.
\]
In particular, $C$ is a rank $N-1$ compression iff
\[
a_1\leq b_1\leq a_2\leq b_2\leq a_3\leq ...\leq a_{N-1}\leq b_{N-1}\leq a_N,
\]
the classic ``interlacing'' of eigenvalues.
\szs
A much more recent example is provided by the theory of \textbf{higher--rank numerical ranges}.
The striking development of this theory was motivated originally by problems in
quantum information theory. Since the introduction of this concept by Choi, Kribs, and \.Zyczkowski [CK\.Z1,CK\.Z2] only a few years ago,
it has indeed been effectively applied in the area of quantum information (see [CPMS\.Z,KPLRdS,LP,LPS1,MM\.Z], for example).
It has also inspired a remarkable
development of its purely mathematical aspects (see, for example, [CHK\.Z,CGHK,Wo,LS,LPS2,DGHP\.Z]). From this point of view the theory of the
higher--rank numerical ranges may be described as a highly successful analysis of \textbf{scalar} compressions of arbitrary
matrices $M\in \dM_N$. This suggests a more general program: characterize the
\textbf{normal} (diagonal) compressions of $M$.
In what follows we begin to carry out this program, although at present the program in its entirety seems out--of--reach.
\szs
The rank--$k$ numerical range of $M$, usually denoted in the literature by $\Lambda_k(M)$, was defined by
Choi, Kribs, and \.Zyczkowski as the set of those complex $\ld$
such that for some rank--$k$ orthogonal projection $P$ we have
\[
      PMP=\ld P.
\]
In terms of compressions, we see that $\ld\in\Lambda_k(M)$ iff
$\ld I_k$ is a (matrix) compression  of $M$. Thus the following fundamental result of Li and Sze [LS] may be
placed in the same family as the Cauchy interlacing theorem (and, in fact, the interlacing theorem plays a role in
the argument of Li and Sze).
\szs
\textbf{Theorem 2:} Given $M\in\dM_N$, let $\ld_j(\theta)$ be an enumeration of the eigenvalues of the (Hermitian)
\[
\mbox{Re}(e^{i\theta}M)=(e^{i\theta}M+e^{-i\theta}M^*)/2
\]
such that
\[
\ld_1(\theta)\leq\ld_2(\theta)\leq\dots\leq\ld_N(\theta).
\]
For each real $\theta$, let the half--plane $H(M,\theta)$ be defined by
\[
H(M,\theta)=e^{i\theta}\{z:\mbox{Re}(z)\leq\ld_{N-k+1}(-\theta)\}.
\]
Then
\beq{E1}{\Lambda_k(M)=\bigcap\{H(M,\theta):\theta\in[0,2\pi]\}.}
\end{equation}
\szs
Our more general program seeks to describe all \textbf{normal} compressions
of $M$, ie to describe those complex $a_1,\dots,a_k$ such that $\diag(a_1,\dots,a_k)$ is a compression of $M$.
Equivalently, we ask when there exist
orthonormal $$u_1,u_2,\dots,u_k$$ such that  $(Mu_i,u_i)=a_i$ for each $i$ and $(Mu_i,u_j)=0$ whenever
$i\neq j$; in particular, $\Lambda_1(M)$ is nothing but the classical numerical range
\[
     W(M)=\{(Mu,u):\|u\|=1\}
\]
(hence the ``higher--rank numerical range'' terminology).
In this work we usually restrict
our attention to the case where $M$ itself is also normal, although we occasionally comment on cases where either $M$ or
its compression may not be normal.
\szs
Note that for normal $M\in\dM_N(\dC)$ Theorem 2 shows that $\Lambda_k(M)$ can be explicitly described in terms of the eigenvalues
$z_1,\dots,z_N$ of $M$:
\beq{E2}{
\Lambda_k(M)=\bigcap_{\#(J)=N-k+1}\conv\{z_j:j\in J\}.}
\end{equation}
We shall refer to this result, first proposed by Choi, Kribs, and \.Zyczkowski, as the CK\.Z conjecture, although
it is now a theorem. The CK\.Z conjecture played an important role in the development of the theory of
higher--rank numerical ranges. For example, while Li and Sze gave an effective description
of $\Lambda_k(M)$ for non--normal $M$ (Theorem 2), their proof of the CK\.Z conjecture was a key step towards the general
result. Of course, the case $k=1$ of (\ref{E2}) is easy and well--known: for normal $M$, $W(M)=\conv\{z_1,\dots,z_N\}$.
\szs
The following observation is often useful.\\
\textbf{Proposition 3:} For every $M\in\dM_N$, if $k\leq N$, $C$ is a rank-$k$ compression of $M$, and $Q$ is
a compression of rank $N-k+1$, then
\[
      W(C)\cap W(Q)\neq\emptyset.
\]
\textbf{Proof:} Let $S$ and $T$ be the subspaces corresponding to compressions $C$ and $Q$. Since the dimensions add
to more than $N$, $S$ and $T$ must intersect non--trivially; let $u$ be a unit vector in $S\cap T$. Then
\[
    (Mu,u)=(Mu,P_Su)=(P_SMu,u)=(Cu,u)\in W(C),
\]
and similarly $(Mu,u)\in W(Q)$. QED
\szs
Applying this observation to the normal case, we see that part of the CK\.Z conjecture is straightforward.\\
\textbf{Proposition 4:} If $M\in\dM_N$ is normal with eigenvalues $z_1,\dots,z_N$, and the rank--$k$ compression $C$
is normal with eigenvalues $c_1,\dots,c_k$, then for every index set $J$ having $\#(J)=N-k+1$
\[
    \conv\{c_1,\dots,c_k\}\cap\conv\{z_j:j\in J\}\neq\emptyset.
\]
In particular,
\[
    \Lambda_k(M)\subseteq \bigcap_{\#(J)=N-k+1} \conv\{z_j:j\in J\}
\]
(compare (\ref{E2})).\\
\textbf{Proof:} We have noted that for normal (finite--dimensional) operators the numerical range is just
the convex hull of the eigenvalues. Thus $W(C)=\conv\{c_1,\dots,c_k\}$.
On the other hand, let $Q$ be the compression to the span of eigenvectors corresponding to $\{z_j:j\in J\}$; then
$Q$ is normal and $W(Q)=\conv\{z_j:j\in J\}$. Apply Proposition 3. In particular, for points $\ld\in\Lambda_k(M)$ we
may let $c_1=c_2=\dots =c_k=\ld$.
QED
\szs
On the other hand, the fact that $\Lambda_k(M)$ completely fills the RHS of (\ref{E2}) is more subtle, in general, although
for certain combinations of $N$ and $k$ it is relatively easy to see. To illustrate this, and to introduce the preoccupations of
the present paper, consider the case $N=5,k=2$. In Figure 1 we see the eigenvalues $z_1,\dots,z_5$ of a normal (in fact, unitary)
$M$ as the outer points of the blue pentagram. It is easy to see that (\ref{E2}) implies that $\Lambda_2(M)$ is the
inner pentagon. As far as we know, there is no simple proof that $\Lambda_k(M)$ fills this pentagon, but three markedly disparate
arguments may be found in the literature:
\szs
(1) in [CHK\.Z] there is an argument based in part on topological concepts such as simple connectivity and winding number;
\szs
(2) as it is easy to conclude (see section 2) that the vertices of the inner pentagon are in $\Lambda_2(M)$, the fact that (whether or not $M$
is normal) $\Lambda_k(M)$ is convex (see [CGHK] and [Wo])) -- a striking extension of the classical Toeplitz--Hausdorff
Theorem for $W(M)$ -- may be used;
\szs
(3) as we have noted, (\ref{E2}) is a direct consequence of the Li and Sze result Theorem 2.
\begin{figure} [htbp]
       \begin{center} \
     \includegraphics[width=10cm,angle=0]{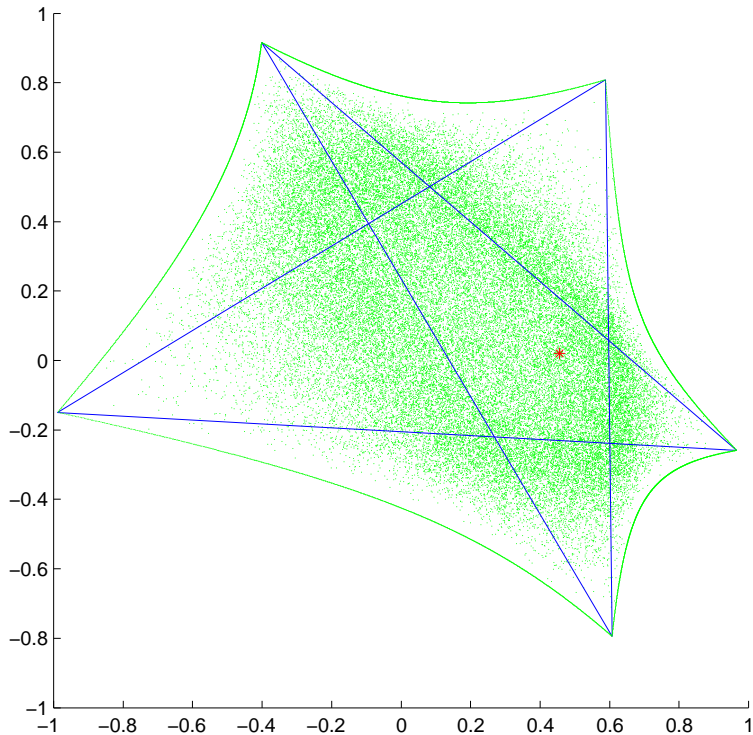}
\caption{Choosing $a$ (red asterisk) at random in $\Lambda_2(M)$  (the inner pentagon), we see that $B(a)$ includes
a ``starfish'' that covers $\Lambda_2(M)$ and more.} \label{fig1}
\end{center}
\end{figure}
\szs
A fourth, and quite different yet again, approach can be obtained by considering those eigenvalue pairs $a,b$ that can belong
to rank--2 normal compressions of $M$. Given $a\in \dC$ we denote by $B(a)$ the set of $b$ that match $a$ in this sense. We shall
prove in section 3 that for $a$ in the inner pentagon $B(a)$ includes a ``starfish'' (outlined in green for the example of Figure 1)
covering the (filled) pentagon (our conjecture, in addition, is that the starfish is precisely $B(a)$). Since $a\in B(a)$ says
that $a\in\Lambda_2(M)$, we conclude once again that $\Lambda_2(M)$ fills the pentagon.
\szs
Plan of the paper: section 2 has some general results, section 3 treats the case $k=2$, section 4 examines continuity of $B(\cdot)$,
and section 5 discusses non--normal compressions.
\szs
Acknowledgements: We have enjoyed many stimulating discussions of matrix compression, particularly
those with M.--D. Choi, C.--K. Li, Y.--T. Poon, N.--S. Sze, and J. F. Queir\'o.
Versions of the material in this paper were developed in [M]. The work of Holbrook and Pereira was supported
in part by Discovery Grants from NSERC of Canada.
\szs
\textbf{2. Some general results (arbitrary $k,N$)}
\szs
Note that if $C$ is a rank--$k$ compression of $M\in\dM_N$ and $C'$ is a rank--$k'$ compression of $C$, then $C'$ is
a rank--$k'$ compression of $M$. Thus Proposition 3 has the following consequence.\\
\textbf{Proposition 5:} If $C$ is a compression of $M\in\dM_N$ then $$W(C)\subseteq W(M).$$
\textbf{Proof:} Regard $z\in W(C)$ as a rank--1 compression $C'$ of $C$, hence of $M$ and apply Proposition 3 with
$k=1$, $C$ replaced by $C'$ and $Q=M$. QED
\szs
Whereas Proposition 4 supplies a \textbf{necessary} condition on the eigenvalues $c_1,\dots,c_k$ of a normal compression
$C$ of normal $M$, the following proposition points out a \textbf{sufficient} condition that is sometimes useful.
An interesting analysis of such necessary vs sufficient conditions may be found in [QD].\\
\textbf{Proposition 6:} If $M\in\dM_N$ is normal with eigenvalues $z_1,\dots,z_N$ then $c_1,\dots,c_k\in\dC$ are
eigenvalues of a normal compression $C$ of $M$ provided that there exists a \textbf{partition} $J_1,\dots,J_k$ of
$\{1,2,\dots,N\}$ such that for each $i=1,\dots,k$
\[
c_i\in\conv\{z_j:j\in J_i\}.
\]
\textbf{Proof:} For each $i$ let $c_i=\sum_{j\in J_i}t_{ij}z_j$ represent $c_i$ as a convex combination. Let $u_1,\dots,u_N$ be an
orthonormal basis of eigenvectors for $M$, with $$Mu_j=z_ju_j.$$ For each $i$, let
\[
w_i=\sum_{j\in J_i} \sqrt{t_{ij}}u_j.
\]
It is easy to check that $w_1,\dots,w_k$ are orthonormal , that $(Mw_i,w_i)=c_i$, and that $(Mw_i,w_h)=0$ if $h\neq i$. It
follows that $C=\diag\{c_1,\dots,c_k\}$ represents the compression of $M$ to the subspace $S=\span\{w_1,\dots,w_k\}$,
ie
\[
C=P_S M|_S.
\]
QED
\szs
In [CK\.Z1] Choi, Kribs, and \.Zyczkowski identified explicitly the higher--rank numerical ranges of Hermitian
matrices, and their argument may be viewed, along the lines of the proof of our next proposition, as an illustration
of the combined force of the necessary condition from Proposition 4 with the sufficient condition from Proposition 6.
Note that the result might also have been
obtained as a special case of the Fan--Pall result, Theorem 1 (taking $b_1=b_2=\dots=b_k$).\\
\textbf{Proposition 7:} If $M\in\dM_N$ is Hermitian with (real) eigenvalues
\[
a_1\leq a_2\leq\dots\leq a_N,
\]
 then for each $k\leq N/2$ we have
\[
\Lambda_k(M) = [a_k,a_{N-k+1}].
\]
If $a_{N-k+1}<a_k$, then $\Lambda_K(M)=\emptyset$.\\
\textbf{Proof:} If $\ld\in\Lambda_k(M)$ then taking $c_1=...=c_k=\ld$ in Proposition 4 we see that
\[
\ld\in\conv\{a_k,\dots,a_N\}=[a_k,a_N].
\]
Likewise, $\ld\in[a_1,a_{N-k+1}]$, so that $\Lambda_k(M)\subseteq [a_k,a_{N-k+1}]$.
\szs
On the other hand, considering the partition of $\{1,\dots,N\}$ into
\[
J_1=\{1,N\},J_2=\{2,N-1\},\dots,J_k=\{k,N-k+1\}
\]
we conclude from Proposition 6 that each $\ld\in[a_k,a_{N-k+1}]$ is in $\Lambda_k(M)$. QED
\szs
As another example of such general arguments we treat the normal compression problem for the case $k=N-1$. This
result goes back to Fan-Pall [FP]; their proof is algebraic in character whereas ours is more geometric. We
restrict to the case where the matrix and its compression have no common eigenvalues since this is where our
general principles are most pertinent; Fan and Pall also treat the general case by means of a direct sum construction.\\
\textbf{Proposition 8:} Let $z_1,\dots,z_N$ and $c_1, \dots,c_{N-1}$ be two collections of complex numbers
having no elements in common. Then there is a normal $M\in\dM_N$ with eigenvalues $z_j$ having a rank--$(N-1)$ normal
compression $C$  with eigenvalues $c_j$ iff the $z_j$ are collinear and alternate with the $c_j$ (in some order)
along the common line.\\
\textbf{Proof:} Let us first show that if such $M,C$ exist then the $z_j$ must be collinear. Label the
$z_j$ lying on the boundary of $W(M)$ in counterclockwise order: $z_1,\dots,z_p$. If the $z_j$ are not
collinear there must be some $z_{k-1},z_k,z_{k+1}$ that are not collinear, as in Figure 2. Proposition 4
requires that $[z_{k-1},z_k]$ meets $W(C)$ at some $\ld$ closest to $z_k$; this $\ld$ is extreme in $W(C)$
and so must be an eigenvalue of $C$. Similarly we have an eigenvalue $\mu$ of $C$ in $[z_k,z_{k+1}]$, as in
Figure 2. Note that Proposition 4 also tells us that $z_k$ cannot be a repeated eigenvalue of $M$, since it
would then coincide with an eigenvalue of $C$.
\szs
Let $u_1,\dots,u_N$ be an orthonormal set of eigenvectors of $M$, with $Mu_j=z_ju_j$, and let orthonormal $v,w$
be eigenvectors of $C$ with $Cv=\ld v$ and $Cw=\mu w$. Expand $v,w$ in terms of the $u_j$:
\[
v=\sum_{j=1}^N a_ju_j,\quad w=\sum_{j=1}^N b_ju_j;
\]
then
\[
\ld=(Cv,v)=(Mv,v)=\sum_{j=1}^N|a_j|^2z_j,
\]
so that $a_j=0$ unless $z_j$ lies on the line through $z_{k-1},z_k$. Similarly $b_j=0$ unless $z_j$ lies on the line
through $z_k,z_{k+1}$. Since $z_k$ is the only common point,
\[
0=(v,w)=a_k\overline{b_k}.
\]
If $a_k=0$ we have $\ld=z_{k-1}$, which we have ruled out, while if $b_k=0$ we have $\mu=z_{k+1}$, also ruled out.

\begin{figure} [htbp]
       \begin{center} \
     \includegraphics[width=10cm,angle=0]{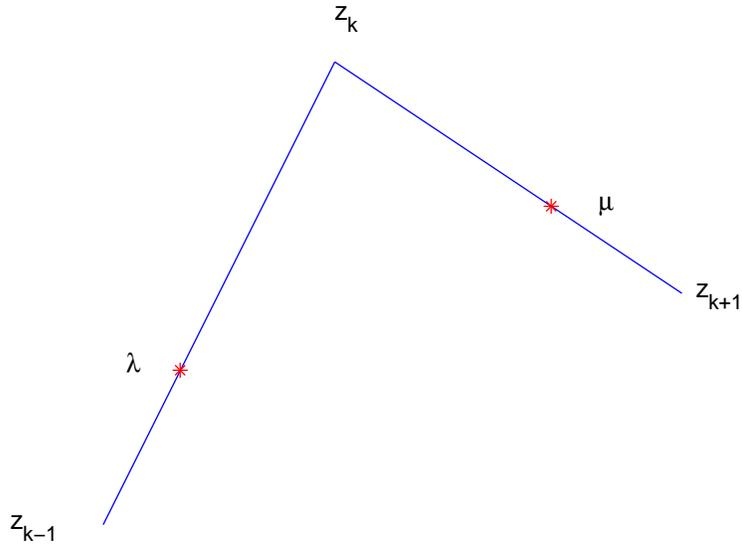}
\caption{An example of the eigenvalue geometry ruled out in the proof of Proposition 8.} \label{fig2}
\end{center}
\end{figure}

\noindent
Thus the eigenvalues all lie on a common line and by an affine map $M\to\alpha I_N + \beta M$ this common line
can be $\dR$, ie we are in the Hermitian case. Proposition 1 then completes the argument, giving the interlacing
property.
\szs
On the other hand, if the collinearity and interlacing conditions are met, the same sort of affine map and
Proposition 1 establish the existence of $M$ and $C$. QED
\newpage\noindent
\textbf{3: Results for $k=2$ and small $N$}
\szs
For $2\times2$ normal compressions $\diag(a,b)$, we can give a more detailed account of the $ab$--geometry, leading up to
an understanding of the ``starfish'' seen in Figure 1.
\szs
Recall that, given normal $M\in\dM_N$ and complex $a$, we denote by $B(a)$ the set of complex $b$ such that $\diag(a,b)$ is
a compression of $M$. Of course, in order that $B(a)$ should be nonempty we must have
\[
a\in\conv\{z_1,z_2,\dots,z_N\},
\]
where the $z_j$ are the eigenvalues of $M$. Note that Proposition 4 also requires that for $b\in B(a)$ we require that the line segment
$[a,b]$ intersect $$\conv\{z_j:j\neq i\}$$ for each $i=1,\dots,N$.
\szs
The simplest case to consider: $N=3$ and the eigenvalues of $M$ form a nontrivial triangle.\\
\textbf{Proposition 9:} Suppose that the eigenvalues $z_1,z_2,z_3$ of normal $M\in\dM_3$ are not collinear. Then $b\in B(a)$ iff
either $a$ is one of these eigenvalues, say $a=z_1$ and $b\in[z_2,z_3]$ (the opposite side of the triangle formed
by $z_1,z_2,z_3$) or $a$ is in one of the sides, say $[z_2,z_3]$, and $b=z_1$.\\
\textbf{Proof:} Since $[a,b]$ must meet each of the triangle's sides, the necessity of the condition is clear. On the
other hand, Proposition 6 shows that these conditions suffice for $a,b$ to be the eigenvalues of a normal compression. QED\\
\textbf{Remark:} Here we have a very simple case of the result of Fan and Pall [FP] where they characterize in general the
case $k=N-1$.
\szs
When $N=4$ we encounter more complex behaviour, such as that seen in Figure 3, where $B(a)$ is a curve interior to
$\conv\{z_1,z_2,z_3,z_4\}$ (except for endpoints).
\newpage
\begin{figure} [htbp]
       \begin{center} \
     \includegraphics[width=10cm,angle=0]{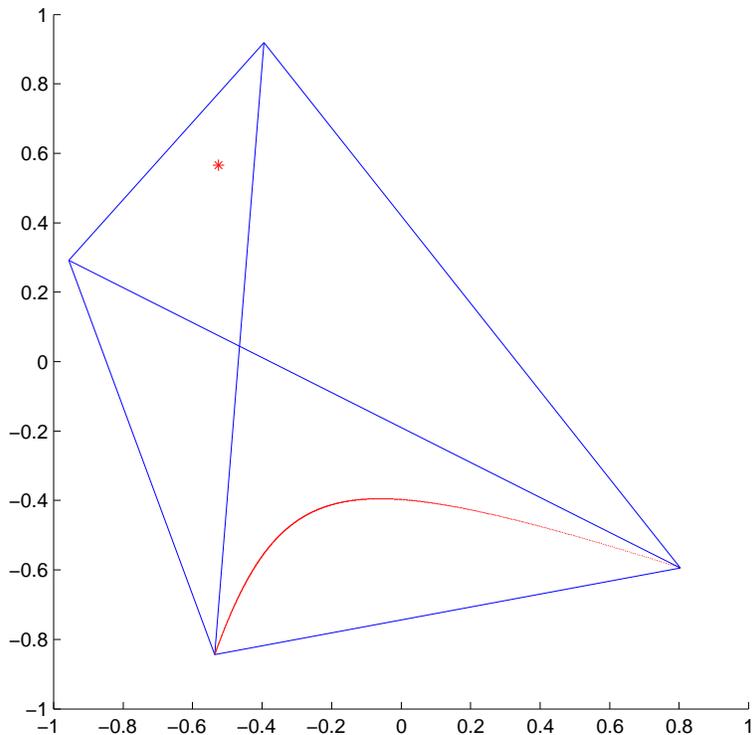}
\caption{For $a$ (red asterisk) strictly inside the upper quadrant (case (a)), we see that $B(a)$ is a curve in
the opposite quadrant.} \label{fig3}
\end{center}
\end{figure}
\noindent
To analyse such behaviour, it will be convenient to assume in what follows that the eigenvalues of $M$ are generic in
the sense that no three are collinear. We may also assume that $M=\diag(z_1,\dots,z_N)$, so that the eigenvectors of
$M$ are the standard basis vectors $e_j$.
\szs
Note that if $b\in B(a)$ we have orthonormal $u,w$ such that
\[
(Mu,u)=a,(Mw,w)=b,\mbox{  and   }(Mu,w)=(Mw,u)=0.
\]
Thus $a=\sum_1^N|u_j|^2z_j$, a convex combination. Let $\Delta_N$ denote the $N$--dimensional simplex, ie
$\conv\{e_1,\dots,e_N\}$; then $|u|^2$ (where the operations are performed componentwise) belongs to
\[
C(a)=\{t\in\Delta_N:a=\sum_1^Nt_jz_j\}.
\]
By exchanging complex arguments between the components of $u$ and $w$ we may assume that $u\geq0$; then the
possible $u$ lie in $\{\sqrt{t}:t\in C(a)\}$. The conditions on $w\in\dC^N$ are then given by
\[
\|w\|=1,w\perp u, w\perp z\circ u, \mbox{ and } w\perp \overline{z}\circ u,
\]
where $\circ$ indicates
Schur (componentwise) multiplication, so that
\[
z\circ u = (z_1u_1,\dots,z_Nu_N)',
\]
with $'$ indicating transpose.
\szs
We may thus describe $B(a)$ as follows.\\
\textbf{Proposition 10:} Given $a\in W(M)(=\conv\{z_1,\dots,z_N\})$,
\[
B(a)=\bigcup_{t\in C(a)} B(a,t),
\]
where
\[
B(a,t)=\{\sum_1^N|w_j|^2z_j:\|w\|=1,w\perp \sqrt{t},z\circ\sqrt{t},\overline{z}\circ\sqrt{t}\}.
\]
\textbf{Proof:} To the discussion above we need only add the observation that
\[
b=(Mw,w)=\sum_1^N|w_j|^2z_j.
\]
QED
\szs
Clearly $C(a)$ is a compact convex subset of $\Delta_N$. It is therefore the convex hull of its extreme points,
which are identified in the following result.\\
\textbf{Proposition 11:}The extreme points of $C(a)$ are those $t\in C(a)$ such that at most three $t_k>0$.\\
\textbf{Proof:} Consider $t\in C(a)$ such that $t_k>0$ for at least four values of $k$. We show that $t$ is \emph{not} extreme.
For convenience assume $t_1,t_2,t_3,t_4>0$. The space
\[
X=\{x\in \dR^N:x_k=0\mbox{ for }k>4\}
\]
is 4--dimensional. Hence
\[
Y=\{x\in X:\sum_1^4x_k=0,\sum_1^4x_k\mbox{Re}(z_k)=0,\sum_1^4x_k\mbox{Im}(z_k)=0\}\neq\{\vec0\}.
\]
Let $\vec0\neq y\in Y$. Then for sufficiently small
$\epsilon>0$ we have $t\pm\epsilon y\in\Delta_N$ and
\[
\sum_k(t\pm\epsilon y)_kz_k=\sum_kt_kz_k=a,
\]
so that $t\pm\epsilon y\in C(a)$. Hence $t$ is not extreme.
\szs
On the other hand, if at most three components, say $t_1,t_2,t_3$ of $t\in C(a)$ are positive, and $t$ is
the average of $t',t''\in C(a)$, then $t'_k,t''_k=0$ for $k>3$. Because no three $z_j$ are collinear,
\[
a=t_1z_1+t_2z_2+t_3z_3
\]
is the \textbf{unique} representation of $a$ as a convex combination of $z_1,z_2,z_3$. Hence $t'=t''=t$.
QED
\szs
For distinct indices $i,j,l$, let $t(i,j,l)$ denote the element of $C(a)$ (if it exists) such that $t_k(i,j,l)=0$
whenever $k\neq i,j,l$. Note that such elements are uniquely determined since
\[
a=t_i(i,j,l)z_i+t_j(i,j,l)z_j+t_l(i,j,l)z_l
\]
represents $a$ uniquely as a point in the triangle $\conv\{z_i,z_j,z_l\}$; here again we use the assumption
that no three of the eigenvalues $z_j$ are collinear. Thus
\beq{E3}{
C(a)=\conv\{t(i,j,l):i,j,l\mbox{ are distinct and }a\in\conv\{z_i,z_j,z_l\}\}.}
\end{equation}
\szs
The complexity of $B(a,t)$ increases with the number of nonzero $t_k$. For example, if only one $t_k>0$,
then $t_k=1$ and $a=z_k$. Here the simple sufficient condition of Proposition 6 is also necessary:
\[
B(a,t)=\conv\{z_j:j\neq k\}.
\]
We see this as follows. Evidently, with $u=\sqrt{t}=e_k$, $u,w$ are orthonormal exactly when $w=\sum_{j\neq k}\alpha_je_j$ with
$\sum_{j\neq k}|\alpha_j|^2=1$; then
\[
b=(Nw,w)=\sum_{j\neq k}|\alpha_j|^2z_j\in\conv\{z_j:j\neq k\},
\]
and any $b\in\conv\{z_j:j\neq k\}$ can be obtained in this way.
\szs
The same sort of simplification occurs if only two or three $t_k>0$.\\
\textbf{Proposition 12:} (a) If $t\in C(a)$ has exactly two positive components, say $t_1,t_2>0$, then
\[
B(a,t)=\conv\{z_j:j>2\}.
\]
(b) If $t\in C(a)$ has exactly three positive components, say $t_1,t_2,t_3>0$, then
\[
B(a,t)=\conv\{z_j:j>3\}.
\]
\textbf{Proof:} (a) Since $a\in\conv\{z_1,z_2\}$, Proposition 6 tells us that
\[
B(a,t)\supseteq\conv\{z_j:j>2\}.
\]
On the other hand, with $u=\sqrt{t}=(\sqrt{t_1},\sqrt{t_2},0,\dots)'$ we see that $u,w$ are orthonormal iff
$\|w\|=1$ and $(w_1,w_2)\perp(\sqrt{t_1},\sqrt{t_2})$; similarly $(Mu,w)=0$ only if
$(w_1,w_2)\perp(\sqrt{t_1}z_1,\sqrt{t_2}z_2)$. Since $z_1\neq z_2$, we have $w_1=w_2=0$ so that
\[
b=(Mw,w)\in\conv\{z_j:j>2\}.
\]
(b) Since $a\in\conv\{z_1,z_2,z_3\}$, Proposition 6 tells us that $$B(a,t)\supseteq\conv\{z_j:j>3\}.$$ On the other hand,
with $u=\sqrt{t}$ we have $u,w$ orthonormal iff $\|w\|=1$ and $$(w_1,w_2,w_3)\perp(\sqrt{t_1},\sqrt{t_2},\sqrt{t_3})$$ and
$(Mu,w)=(Mw,u)=0$ only if
\[
(w_1,w_2,w_3)\perp(\sqrt{t_1}\mbox{Re}(z_1),\sqrt{t_2}\mbox{Re}(z_2),\sqrt{t_3}\mbox{Re}(z_3)),
(\sqrt{t_1}\mbox{Im}(z_1),\sqrt{t_2}\mbox{Im}(z_2),\sqrt{t_3}\mbox{Im}(z_3)).
\]
Since $z_1,z_2,z_3$ are not collinear,
\[
(1,1,1),\quad(\mbox{Re}(z_1),\mbox{Re}(z_2),\mbox{Re}(z_3)),\quad(\mbox{Im}(z_1),\mbox{Im}(z_2),\mbox{Im}(z_3))
\]
are linearly independent.
We must have $w_1=w_2=w_3=0$ so that $b=(Mw,w)\in\conv\{z_j:j>3\}$. QED
\szs
We are now in a position to understand the features of Figure 3 and, indeed, to analyse all the possibilities
when $N=4$. We treat in detail the case where $z_1,z_2,z_3,z_4$ are all extreme in $\conv\{z_1,z_2,z_3,z_4\}$;
the case where one of the eigenvalues lies in the interior of $W(M)$ (eg $z_4\in\conv\{z_1,z_2,z_3\}$) can be
treated similarly.\\
\textbf{Proposition 13:} Let $N=4$ and suppose that $z_1,z_2,z_3,z_4$ are all extreme in $W(M)$ and are
numbered in counterclockwise order. The diagonals $[z_1,z_3]$ and $[z_2,z_4]$ meet at $q$ and divide $W(M)$ into
four quadrants. Consider $a\in W(M)$; the possibilities for $B(a)$ are as follows.\\
(a) See figure 3: $a$ lies in the interior of one of the quadrants. For convenience, assume that $a\in\conv\{z_1,z_2,q\}$;
let $x=t(1,2,3)$, $y=t(1,2,4)$. Then $B(a)$ is the curve traced out by the function $b(r)$ defined for $0<r<1$ by
\[
b(r)=\sum_{k=1}^4\frac{(x_k-y_k)^2}{(1-r)x_k+ry_k}z_k\Big/ \sum_{k=1}^4\frac{(x_k-y_k)^2}{(1-r)x_k+ry_k}.
\]
Note that $x_4=0$ and $y_3=0$ so that
\[
\lim_{r\to0}b(r)=z_4, \quad \lim_{r\to1}b(r)=z_3,
\]
and we obtain a continuous curve parametrized on $[0,1]$ when we interpret $b(0)$ as $z_4$ and $b(1)$ as $z_3$. Except for
these endpoints, the curve lies in the interior of the opposite quadrant $\conv\{z_3,z_4,q\}$.\\
(b) If $a$ lies in the interior of one of the sides of $W(M)$ then $B(a)$ is the opposite side (eg if $a$ is inside
$[z_1,z_2]$ then $B(a)=[z_3,z_4]$). If $a=z_k$ then $B(a)$ is the opposite triangle $\conv\{z_j:j\neq k\}$.\\
(c) See Figure 4: $a$ lies interior to the diagonals but is not $q$; say $a$ is interior to $[z_1,q]$. Then $B(a)$
is the T--shaped object $[z_2,z_4]\cup[q,z_3]$.\\
(d) If $a=q$ then $B(a)$ is the union of the two diagonals.
\newpage

\begin{figure} [htbp]
       \begin{center} \
     \includegraphics[width=10cm,angle=0]{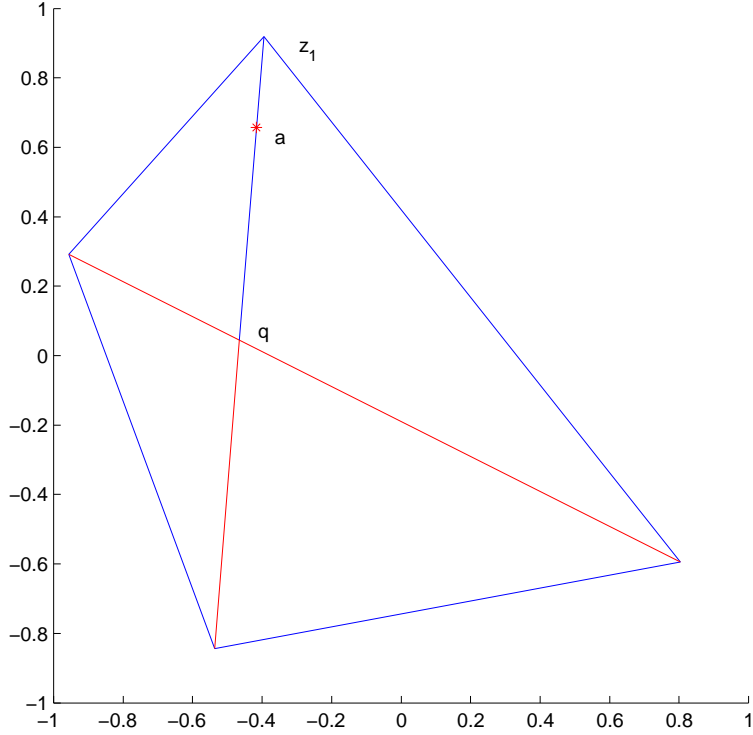}
\caption{For $a$ (red asterisk) strictly inside the segment $[z_1,q]$ (case (c)), we see that $B(a)$ is the T--shaped
object consisting of $[z_2,z_4]\cup[q,z_3]$.} \label{fig4}
\end{center}
\end{figure}
\noindent
\textbf{Proof:} (a) Since $a$ lies in the triangles $\conv\{z_1,z_2,z_3\}$ and $\conv\{z_1,z_2,z_4\}$ but in no other
triangle of
eigenvalues, $C(a)=[t(1,2,3),t(1,2,4)]=[x,y]$ (recall the relation (\ref{E3})). For $0<r<1$ consider the $t\in C(a)$ given by
$t=(1-r)x+ry$. We shall see that $B(a,t)$ consists of the single point $b(r)$. We take
$u=\sqrt{t}$ and note that the conditions on $w$ are: $w\perp u$, $w\perp u\circ \mbox{Re}(z)$,
$w\perp u\circ \mbox{Im}(z)$, and $\|w\|=1$. Thus $w\circ\sqrt{t}\perp\vec1_4,\mbox{Re}(z),\mbox{Im}(z)$,
where $\vec1_4$ denotes $[1,1,1,1]$. Again we invoke
linear independence of $\vec1_4,\mbox{Re}(z),\mbox{Im}(z)$: $w\circ\sqrt{t}$ lies in the one--dimensional space
\[
dC^4\circleddash\span\{\vec1_4,\mbox{Re}(z),\mbox{Im}(z)\}.
\]
There is a natural choice of (nonzero) vector in this space: $x-y$ (because $(x,\vec1_4)=(y,\vec1_4)=1$,
$(x,\mbox{Re}(z))=(y,\mbox{Re}(z))=\mbox{Re}(a)$, and $(x,\mbox{Im}(z))=(y,\mbox{Im}(z))=\mbox{Im}(a)$). Thus
\[
w=\alpha(x-y)/\circ\sqrt{t},
\]
where $/\circ$ indicates entrywise division and $\alpha$ is some complex number. Recalling that $\|w\|=1$, we derive our
formula for $(Mw,w)=b(r)$.
\szs
The necessary condition of Proposition 4 shows that the curve (ie $B(a)$) lies in both $\conv\{z_2,z_3,z_4\}$ and $\conv\{z_1,z_3,z_4\}$, so
that it must lie in the (closed) opposite quadrant $\conv\{z_3,z_4,q\}$. To see that the curve (except for endpoints) lies in the
interior of that quadrant, examine the arguments below, showing that for $b$ on the quadrant boundary (except for $z_3$ and $z_4$)
a matching $a$ cannot be interior to the upper quadrant, and note that $b\in B(a)$ iff $a\in B(b)$.
\szs
(b) If $a$ is interior to one of the sides, say $[z_1,z_2]$, then $C(a)$ consists of a single $t$ with two positive components; apply
Proposition 12(a) to see that $B(a)=B(a,t)=[z_3,z_4]$. If $a=z_1$, Propositions 4 and 6 imply that $B(a)=\conv\{z_2,z_3,z_4\}$.
\szs
(c) Suppose $a$ is interior to $[z_1,q]$; then the relation (\ref{E3}) tells us that
\[
C(a)=\conv\{t(1,2,3),t(1,3,4),t(1,2,4)\}.
\]
Let $t(1,2,3)=x=[x_1,0,x_3,0]'$; this is also $t(1,3,4)$. Let $t(1,2,4)=y=[y_1,y_2,0,y_4]'$, so that
$C(a)=\{t(r):0\leq r\leq1\}$, where
\[
t(r)=[(1-r)x_1+ry_1,ry_2,(1-r)x_3,ry_4]'.
\]
For $r=0$, Proposition 12(a) tells us that $B(a,t(0))=[z_2,z_4]$, while for $0<r\leq1$ we claim that
$B(a,t(r))$ is a single point $b(r)$ that moves along $[z_3,q)$, covering it completely. Indeed, reasoning as in
(a), we see that $b(r)=(Mw(r),w(r))$ where $w(r)$ is a normalized version of
\[
(t(0)-t(1))/\circ\sqrt{t(r)}.
\]
Note that $w_2(r),w_4(r)$ are proportional to $-y_2/\sqrt{ry_2},-y_4/\sqrt{ry_4}$ respectively, so that
\[
\frac{|w_2(r)|^2}{|w_4(r)|^2}=\frac{y_2}{y_4}.
\]
Since $a=y_1z_1+y_2z_2+y_4z_4$ lies on $[z_1,z_3]$, we conclude that $b(r)\in[z_1,z_3]$ also. The necessary
condition of Proposition 4 then tells us that $b(r)\in[z_3,q]$. Since $t_2(r)$ and $t_4(r)$ tend to 0 as
$r\to0$, $\lim_{r\to0}b(r)=q$. Moreover, $t(1)=[y_1,y_2,0,y_4]$ so that Proposition 12(b) implies that
$b(1)=z_3$. Finally, since $b(r)$ is continuous over $0<r\leq1$, its values cover $[z_3,q)$.
\szs
(d) This case may be treated by an argument rather similar to that of (c). QED
\szs
We now have the tools to continue the theme of Proposition 12, treating the case when exactly \textbf{four}
of the components of $t\in C(a)$ are positive.\\
\textbf{Proposition 14:} Suppose that $N>4$ and that $t\in C(a)$ has exactly four positive components; for
convenience, assume that $t_1,t_2,t_3,t_4>0$ and that $a$ lies in the upper quadrant relative to
$Q=\conv\{z_1,z_2,z_3,z_4\}$, ie $a$ is interior to $\conv\{z_1,z_2,q\}$ (see Figure 3, with the understanding
that it is now intended to show only the relation of $a$ to $z_1,z_2,z_3,z_4$, and Proposition 13). Let $\beta$
be the curve traced out by $b(\cdot)$ of Proposition 13(a) (and shown in Figure 3). Then
\[
B(a,t)=\conv\{\beta,z_5,z_6,\dots,z_N\}.
\]
\textbf{Proof:} With $u=\sqrt{t}$, we see that the conditions on $w$, namely
\[
w\perp u,u\circ\mbox{Re}(z),u\circ\mbox{Im}(z)\mbox{  and  }\|w\|=1,
\]
reduce to
\[
\tilde{w}\perp\tilde{u},\tilde{u}\tilde{\circ\mbox{Re}(z)},\tilde{u}\tilde{\circ\mbox{Im}(z)},
\]
where $\tilde{w}=(w_1,w_2,w_3,w_4)'$, $\tilde{u}=(u_1,u_2,u_3,u_4)'$ etc,
and
\[
\|\tilde{w}\|^2+\sum_{k>4}|w_k|^2=1.
\]
Thus $\tilde{w}/\|\tilde{w}\|$ is subject to the same conditions as $w$ in the proof of Proposition 13(a).
It follows that
\[
(Mw,w)=\|\tilde{w}\|^2b(r)+\sum_{k>4}|w_k|^2z_k
\]
where $b(r)$ can be any point on the curve $\beta$. QED
\szs
Proposition 14 allows us to understand, in large part, the phenomenon illustrated in Figure 1. Let $N=5$ and
suppose that each eigenvalue $z_k$ is an extreme point of $W(M)=\conv\{z_1,\dots,z_5\}$ (eg whenever $M$ is unitary).
For convenience, label the $z_k$ in counterclockwise order. Suppose that $a$ lies strictly inside the central pentagon
(which is known to be $\Lambda_2(M)$ in this case). For each $k$ let $\beta_k$ denote the curve obtained as in Proposition 14 by regarding
$a$ as an element of the quadrilateral $Q_k=\conv\{z_j:j\neq k\}$. Note that $\beta_k$ connects $z_{k+2}$ and $z_{k+3}$
(numbering modulo 5) and lies in the quadrant of $Q_k$ opposite to the one containing $a$. We claim that (as illustrated
in Figure 1) $B(a)$ includes the whole ``starfish'' region bounded by $\beta_1,\beta_2,\dots,\beta_5$.
\szs
To see this note that the starfish is the union of the wedges $W_k=\conv\{\beta_k,z_k\}$, so it suffices to show that each $W_k\subseteq B(a)$.
Since $a\in Q_k$ there is $t\in C(a)$ such that $t_k=0$. Then Proposition 14 tells us that $B(a,c)=W_k$.
\szs
Figure 1 was obtained by first computing $C(a)$ via the relation (\ref{E3}) as
$$\conv\{t(k,k+2,k+3):k=1,2,\dots,5\}$$
(note that for $a$ in the inner pentagon, the only eigenvalue triangles containing $a$ correspond to the triples
$z_k,z_{k+2},z_{k+3}$). To generate each of the thousands of $b$'s in $B(a)$, plotted as green points in Figure 1,
our MATLAB program first chose a ``random'' point $t\in C(a)$ (ie a random convex combination of the five
$c(k,k+2,k+3)$), put $u=\sqrt{t}$, then computed $b=(Nw,w)$ where $w$ was chosen ``randomly'' in
$$\dC^5\circleddash\span\{u,u\circ\mbox{Re}(z),u\circ\mbox{Im}(z)\}$$
(and normalized so that $\|w\|=1$). The curves $\beta_k$ were added using the formula of Proposition 13(a).
Such simulations strongly suggest the following ``starfish conjecture'', since no green dots fall outside the starfish: in such
a situation (and in particular when $N=5$ and $M$ is unitary), $B(a)$ not only contains the starfish but is equal to it.
\szs
We have seen in the discussion of Figure 1 that for $N=5$ and $a,b\in\Lambda_2(M)$ we always have $a,b$ as eigenvalues of
a normal compression of $M$. The following proposition points out that this is true for any $N$ -- and that $N=5$ is, in fact,
the only subtle case.\\
\textbf{Proposition 15:} Let $M$ be normal in $\dM_N$ and such that the eigenvalues $z_1,\dots,z_N$ are distinct and
each is an extreme point of $W(M)$ (eg $M$ unitary). Then $a,b\in\Lambda_2(M)$ implies that $\begin{bmatrix}a&0\\0&b\end{bmatrix}$
is a compression of $M$.\\
\textbf{Proof:} For $N\leq3$, $\Lambda_2(M)=\emptyset$. For even $N\geq4$, the relation (\ref{E2}) tells us that $\Lambda_2(M)$
is the ``inner $N$--gon'' cut off by the line segments $[z_j,z_{j+2}]$ (indexing modulo $N$). Thus for even $N\geq4$
\[
\Lambda_2(M)=\conv\{z_j: j\mbox{ odd}\}\cap\conv\{z_j: j\mbox{ even}\},
\]
and Proposition 6 suffices. For $N=5$ the ``starfish'' discussion
proves our assertion. For odd $N\geq7$ we see that $\conv\{z_j:
j\mbox{ odd}\}\supseteq\Lambda_2(M)$ and $\conv\{z_j: j\mbox{
even}\}$ covers all of $\Lambda_2(M)$ except that part lying in
$Q=\conv\{z_1,z_2,z_{N-1},z_N\}$. Hence Proposition 6 suffices for
$a\not\in Q,b\in\Lambda_2(M)$. The same argument applies for
$a\not\in\tilde{Q}=\conv\{z_2.z_3,z_4,z_5\}$ and because $N>5$ this
covers any $a\in Q$. QED
\szs
\textbf{4. Continuity of $B(\cdot)$}
\szs
A natural assertion of ``continuity'' for $B(\cdot)$ might be that $d_H(B(a'),B(a))\to0$ as $a'\to a$, where
$d_H(X,Y)$ is the Hausdorff distance between compact nonempty sets $X,Y\subset\dC$. Recall that
\[
d_H(X,Y)=\max\{\hat{d}_H(X,Y),\hat{d}_H(Y,X)\},
\]
where
\[
\hat{d}_H(X,Y)=\max_{x\in X}(\min_{y\in Y}|x-y|).
\]
However, we have seen simple examples where this fails: recall the analysis of $B(a)$ for various
$a\in\conv\{z_1,z_2,z_3,z_4\}$ that was provided by Proposition 13. If $a'$ lies in the interior of
$[z_1,z_2]$ and $a'\to a=z_1$, then $B(a')=[z_3,z_4]$ ``jumps'' to $B(a)=\conv\{z_2,z_3,z_4\}$.
A perhaps more surprising example: let $a$ be interior to $[z_1,q]$ as in Figure 4; for $a'$ approaching $a$
from the interior of $\conv\{z_1,z_2,q\}$ we see $B(a')$ as a curve joining $z_3$ and $z_4$ in $\conv\{z_3,z_4,q\}$,
whereas for $a'$ approaching $a$
from the interior of $\conv\{z_1,z_4,q\}$ we see $B(a')$ as a curve joining $z_2$ and $z_3$ in $\conv\{z_2,z_3,q\}$.
\szs
In spite of such ``failures'' we'll show that $B(\cdot)$ is continuous with respect to Hausdorff distance at most points
of $W(M)$ and enjoys a ``one--sided'' Hausdorff continuity in general.
\szs
Our standard set--up for this discussion is as in section 3, ie we assume $M$ is normal in $\dM_N$ and is in diagonal
form: $M=\diag(z)$, where no three eigenvalues are collinear. Thus $W(M)=\conv\{z_1,\dots,z_N\}$ and $B(a')=\emptyset$ if
$a'\not\in W(M)$. Seeking continuity, we restrict attention to $a'\to a$ with $a',a\in W(M)$. Note that if $N=3$ and
$a'$ is interior to $W(M)=\conv\{z_1,z_2,z_3\}$, we again have $B(a')=\emptyset$, since $b\in B(a')$ and Proposition 4
would require that $[a',b]$ meet each side of the triangle $W(M)$. We therefore restrict also to cases where $N\geq4$.
\szs
\textbf{Proposition 16:} If $N\geq 4$, $B(a)$ is a compact nonempty set for any $a\in W(M)$.\\
\textbf{Proof:} Let $t\in C(a)$. Since $N\geq4$,
\[
\dC^N\circleddash\span\{\sqrt{t},\sqrt{t}\circ\mbox{Re}(z),\sqrt{t}\circ\mbox{Im}(z)\}
\]
is nontrivial ($\not=\{\vec0\}$). Let $w$ be a unit vector in this space; then $b=(Mw,w)\in B(a,t)$,
so $B(a)\not=\emptyset$.
\szs
For compactness, consider $b_n\in B(a)$; there exist orthonormal pairs $u_n,w_n$ such that
\[
(Mu_n,u_n)=a,\quad (Mw_n,w_n)=b_n,\quad (Mu_n,w_n)=(Mw_n,u_n)=0.
\]
Since the sequences $u_n,w_n$ are bounded, local compactness in $\dC^{2N}$ implies that,
for some subsequence $n_k$,
\[
u_{n_k}\to_k u,\quad w_{n_k}\to_k w.
\]
Then $u,w$ are orthonormal and
\[
(Mu,u)=a,\quad (Mw,w)=\lim_k b_{n_k}=b,\quad (Mu,w)=(Mw,u)=0.
\]
The limit point $b$ is in $B(a)$. QED
\szs
A related argument shows that, in general, $B(\cdot)$ is continuous in a one--sided Hausdorff sense.\\
\textbf{Proposition 17:} If $a,a_n\in W(M)$ and $a_n\to a$, then
\beq{E4}{
\hat{d}_H(B(a_n),B(a))\to_n0.}
\end{equation}
\textbf{Proof:} Recall that $\hat{d}_H(X,Y)=\max_{x\in X}(\min_{y\in Y}|x-y|)$. Thus, if (\ref{E4}) were
to fail we'd have some $\epsilon>0$, subsequence $n_k$, and $b_k\in B(a_{n_k})$ such that for all $b\in B(a)$
\[
|b_k-b|\geq\epsilon.
\]
By restricting to such a subsequence we may assume that $b_n\in B(a_n)$. Let $u_n,w_n$ be orthonormal pairs
such that
\[
(Mu_n,u_n)=a_n,\quad (Mw_n,w_n)=b_n,\quad (Mu_n,w_n)=(Mw_n,u_n)=0.
\]
There is a subsequence $n_k$ such that
\[
u_{n_k}\to_k u,\quad w_{n_k}\to_k w.
\]
Hence $u,w$ are orthonormal and
\[
(Mu,u)=\lim_k a_{n_k}=a,\quad (Mw,w)=\lim_k b_{n_k}=b,\quad (Mu,w)=(Mw,u)=0.
\]
It follows that $b=\lim_k b_{n_k} \in B(a)$, contradicting $|b_{n_k}-b|\geq\epsilon$. QED
\szs
In terms of the obvious extension of Hausdorff distance to compact nonempty subsets of $\Delta_N$, we note
that $C(\cdot)$ is continuous and in fact satisfies a Lipschitz condition for each fixed $M$.\\
\textbf{Proposition 18:} There is a constant $K<\infty$ depending only on $M$ such that for all $a,a'\in W(M)$
\[
d_H(C(a),C(a'))\leq K|a-a'|.
\]
\textbf{Proof:} For each triple $i,j,k$ of distinct indices, we have assumed that $z_i,z_j,z_k$ are not collinear.
Thus the matrix
\[
T=\begin{bmatrix}1&1&1\\\mbox{Re}(z_i)&\mbox{Re}(z_j)&\mbox{Re}(z_k)\\\mbox{Im}(z_i)&\mbox{Im}(z_j)&\mbox{Im}(z_k)
\end{bmatrix}
\]
is nonsingular. Given $a\in\conv\{z_i,z_j,z_k\}$, consider $t_{ijk}=t(i,j,k)$ as in (\ref{E3}).
Let $\hat{t}_{ijk}$ be the vector in $\dR^3$ recording the $i$, $j$, $k$ --components of $t_{ijk}$, ie
the only components that may be positive. We have $T\hat{t}_{ijk}=(1,\mbox{Re}(a),\mbox{Im}(a))'$ so that
$\hat{t}_{ijk}=T^{-1}(1,\mbox{Re}(a),\mbox{Im}(a))'$. In terms of the operator norm $\|T^{-1}\|$ we have
\[
\|\hat{t}_{ijk}-\hat{t}'_{ijk}\|\leq\|T^{-1}\||a-a'|
\]
for any other $a'\in\conv\{z_i,z_j,z_k\}$. Let $K$ be the maximum of $\|T^{-1}\|$ over all such triangles
$\conv\{z_i,z_j,z_k\}$.
\szs
The line segments $[z_i,z_j]$ form a ``grid'' criss--crossing $W(M)$, dividing it into regions. Suppose $a,a'$ lie
in the same one of these regions (boundary points allowed). Then the set $Q$ of triples $i,j,k$ such that
$a\in\conv\{z_i,z_j,z_k\}$ is the same as that for $a'$. In view of (\ref{E3}), each $t\in C(a)$ can be expressed as a
convex combination
\[
t=\sum_{ijk\in Q} s_{ijk}t_{ijk}.
\]
Putting
\[
t'=\sum_{ijk\in Q} s_{ijk}t'_{ijk},
\]
we have $t'\in C(a')$ and $\|t-t'\|\leq K|a-a'|$. The roles of $a,a'$ may be reversed, so we see that if $a,a'$ are
in the same region (boundary points allowed),
\[
d_H(C(a),C(a'))\leq K|a-a'|.
\]
Finally, for any $a,a'\in W(M)$, the line segment $[a,a']$ intersects the grid in a sequence of points
$a_0,a_1,\dots,a_n$ ordered along $[a,a']$ with $a_0=a,a_n=a'$. By the argument above,
\[
d_H(C(a_k),C(a_{k+1}))\leq K|a_k-a_{k+1}|,
\]
so that ($d_H$ is a metric)
\[
d_H(C(a),C(a'))\leq K\sum_{k=0}^{n-1}|a_k-a_{k+1}|=K|a-a'|.
\]
QED
\szs
Next we show that $B(\cdot)$ is $d_H$--continuous at any point that is ``off the grid'', and that continuity
is uniform if we stay bounded away from the grid.\\
\textbf{Proposition 19:} If $a\in W(M)$ but $a$ does not lie on any line segment $[z_i,z_j]$, then $a'\to a$
implies that
\[
d_H(B(a'),B(a))\to 0.
\]
In fact, on any subset $S(d)\subset W(M)$ that is a positive distance $d$ from the grid
\[
G=\bigcup\{[z_i,z_j]:i,j=1,\dots,N\},
\]
so that
\[
S(d)=\{a\in W(M):\min_{g\in G}|a-g|\geq d\},
\]
the map $a\mapsto B(a)$ is uniformly continuous.\\
\textbf{Proof:} In this discussion $i,j,k$ always denotes a triple of distinct indices. Let
\[
Q=\bigcup\{C(a):a\in S(d)\};
\]
we claim that
\[
\min_{t\in Q}(\max_{i,j,k} t_it_jt_k)
\]
is positive. Otherwise, by compactness, we'd have some $a\in S(d)$ and $t\in C(a)$ such that $\max_{i,j,k} t_it_jt_k=0$. This
can only happen if $t$ has at most two positive components, say $t_i,t_j$; then $a\in[z_i,z_j]$, which we have ruled out.
\szs
Given linearly independent $q,r,s\in\dC^N$, let $P(q,r,s)$ denote orthogonal projection onto
\[
\dC^N\circleddash\span\{q,r,s\}.
\]
The map $(q,r,s)\mapsto P(q,r,s)$ is uniformly continuous if we ``stay away from dependence''; to be precise, for any $0<h<H<\infty$
this map is uniformly continuous on
\[
Q(h,H)=\{(q,r,s):\|q\|,\|r\|,\|s\|\leq H, \max_{i,j,k}|\det\begin{bmatrix}q_i&q_j&q_k\\r_i&r_j&r_k\\s_i&s_j&s_k\end{bmatrix}|\geq h\}.
\]
Now the values $(\sqrt{t},\sqrt{t}\circ\mbox{Re}(z),\sqrt{t}\circ\mbox{Im}(z))$ where $t\in Q$ lie in some fixed $Q(h,H)$
because each
\[
\det\begin{bmatrix}1&1&1\\ \mbox{Re}(z_i)&\mbox{Re}(z_j)&\mbox{Re}(z_k)\\\mbox{Im}(z_i)&\mbox{Im}(z_j)&\mbox{Im}(z_k)\end{bmatrix}
\]
is nonzero, so that
\[
\max_{i,j,k}|\sqrt{t_it_jt_k}\det\begin{bmatrix}1&1&1\\ \mbox{Re}(z_i)&\mbox{Re}(z_j)&\mbox{Re}(z_k)\\
\mbox{Im}(z_i)&\mbox{Im}(z_j)&\mbox{Im}(z_k)\end{bmatrix}|\geq h
\]
for some positive $h$. Thus the map $t\mapsto P(\sqrt{t},\sqrt{t}\circ\mbox{Re}(z),\sqrt{t}\circ\mbox{Im}(z))=P[t]$ is
uniformly continuous on $Q$: given $\epsilon_1>0$ there is $\delta_1>0$ such that $t,t'\in Q$ and $\|t-t'\|\leq\delta_1$
implies $\|P[t]-P[t']\|\leq\epsilon_1$.
\szs
In view of Proposition 18, there is $\delta>0$ such that $|a-a'|\leq\delta$ implies $d_H(C(a),C(a'))\leq\delta_1$. Consider
$b\in B(a)$; for some $t\in C(a)$ we have $b\in B(a,t)$ so that $b=(Mw,w)$ for some unit $w$ with $P[t]w=w$. Let $t'\in C(a')$
be such that $\|t-t'\|\leq\delta_1$; then $\|w-P[t']w\|\leq\epsilon_1$. Note that
\[
1-\epsilon_1\leq\|P[t']w\|\leq1,
\]
and let $w'=P[t']w/\|P[t']w\|$; then $b'=(Mw',w')\in B(a')$ and
\[
|b-b'|=|(Mw,w)-(Mw',w')|\leq 2\|M\|\,\|w-w'\|.
\]
It is easy to see that $\|w-w'\|\leq 2\epsilon_1/(1-\epsilon_1)$, so that given any $\epsilon>0$ we have $|b-b'|\leq\epsilon$
by an appropriate choice of $\epsilon_1$. We have shown that $|a-a'|\leq\delta$ implies that $\hat{d}_H(B(a),B(a'))\leq\epsilon$.
Since the roles of $a,a'$ may be reversed, we also have $d_H(B(a),B(a'))\leq\epsilon$. QED
\szs
Note that \textbf{sometimes} $B(\cdot)$ is continuous even at points that are on the grid. For example, from Proposition 13(a) and 13(b)
we can see that there is continuity everywhere on the boundary segments $[z_i,z_{i+1}]$ except at the endpoints.
\szs
\textbf{5. Related results}
\szs
We offer some remarks on the apparently more difficult problem of characterizing \textbf{arbitrary} compressions
of a normal matrix $M$. Suppose again that $M$ is $N\times N$, and is represented by the diagonal matrix $\diag(z)$ and that
$X$ is a rank--$k$ compression of $M$, ie there is a $k$--dimensional subspace $S$ such that $X=P_SM|_S$. From Proposition 3
we obtain a \textbf{necessary} condition on $X$: the (classical) numerical range
$W(X)$ of $X$ must intersect the convex hull of any subset of the eigenvalues $z_j$ having size $N-k+1$.
\szs
When $k=2$, ie $X$ is represented by a $2\times2$ matrix, the numerical range $W(X)$ determines $X$ uniquely as an
operator. Indeed, $W(X)$ is a (filled--in) ellipse in this case with the eigenvalues of $X$ as foci and the length of
the minor axis is the modulus of the off--diagonal entry of any upper--triangular matrix for $X$. Let's consider the
problem of characterizing such compressions $X$ geometrically via the elliptical $W(X)$ in the cases where $N=3$ and $N=4$.
\szs
When $N=3$, the necessary condition of above tells us that $W(X)$ must be tangent to each of the three sides of $\conv\{z_1,z_2,z_3\}$
(recall that Proposition 5 tells us that in general we must have $W(X)\subseteq W(M)=\conv\{z_j:j=1,\dots,n\}$).
In fact, Williams showed long ago that the
necessary condition is also sufficient when $N=3$ (see [Wi]).
\szs
When $N=4$ we consider the case where the eigenvalues $z_j$ form a quadrilateral $Q$. The necessary condition above tells
us that $W(X)$ must intersect each of the four triangles $T_i=\conv\{z_j:j\neq i\}$. Thus $W(X)$ must intersect each of the
quadrants $T_i\cap T_k$. This phenomenon is borne out by numerical experiments such as Figure 5 illustrates, but it is not
clear what additional conditions must be satisfied by $W(X)$, even in this $N=4$ case. Of course, if by chance $W(X)$ is tangent
to all three sides of some $T_i$, then Williams' result tells us that $X$ is indeed a 2--dimensional compression.
\begin{figure} [htbp]
       \begin{center} \
     \includegraphics[width=10cm,angle=0]{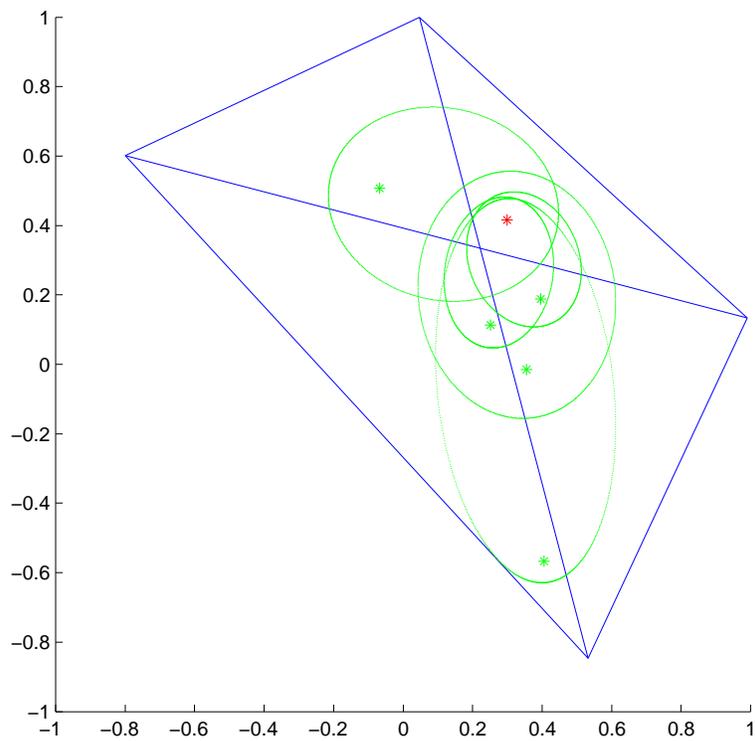}
\caption{Shows the (elliptical) boundaries of the numerical ranges of several (nonnormal) compressions
of a $4\times4$ normal $M$, each  compression having $a$ (red asterisk) as an eigenvalue (therefore  seen as one of the
foci of each ellipse)} \label{fig5}
\end{center}
\end{figure}
\newpage\noindent
References:
\szs
[Cau] A. L. Cauchy, Sur l'\'equation \`a l'aide de laquelle on d\'etermine les in\'egalit\'es
s\'eculaires des mouvements des plan\`etes, Oeuvres compl\`etes, Second Ser., IX, 174--195
\szs
[CK\.Z1] M.--D. Choi, D. W. Kribs, and K. \.Zyczkowski, Higher--rank numerical ranges and compression problems,
Linear Algebra Appl. 418, 828--839, 2006
\szs
[CK\.Z2] M.--D. Choi, D. W. Kribs, and K. \.Zyczkowski, Quantum error correcting codes from the compression
formalism, Rep. Math. Phys. 58, 77--91, 2006
\szs
[CHK\.Z] M.--D. Choi, J. A. Holbrook, D. W. Kribs, and K. \.Zyczkowski, Higher--rank numerical ranges of unitary
and normal matrices, Operators and Matrices 1, 409--426, 2007
\szs
[CGHK] M.--D. Choi, M. Giesinger, J. A. Holbrook, and D. W. Kribs, Geometry of higher--rank numerical ranges,
Linear and Multilinear Algebra 56, 53--64, 2008
\szs
[DGHP\.Z] C. F. Dunkl, P. Gawron, J. Holbrook, Z. Puchala, and K. \.Zyczkowski, Numerical shadows: measures
and densities on the numerical range, Linear Algebra Appl. 434, 2042--2080, 2011
\szs
[FP] K. Fan and G. Pall, Imbedding conditions for Hermitian and normal matrices, Canadian J. Math. 9, 298--304, 1957
\szs
[GPMS\.Z] P. Gawron, Z. Puchala, J. Miszczak, L. Skowronek, and K. \.Zyczkowski, Restricted numerical range:
a versatile tool in the theory of quantum information, J. Math. Physics 51, 2010
\szs
[KPLRdS] D. W. Kribs, A. Pasieka, M. Laforest, C. Ryan, and M. P. da Silva, Research problems on numerical ranges
in quantum computing, Linear and Multilinear Algebra 57, 491-502, 2009
\szs
[LP] C.--K. Li and Y.--T. Poon, Generalized numerical ranges and quantum error correction, J. Operator Theory
66, 335--351, 2011
\szs
[LPS1] C.--K. Li, Y.--T. Poon, and N.--S. Sze, Higher rank numerical ranges and low rank perturbations of
quantum channels, J. Math. Analysis Appl. 348, 843--855, 2008
\szs
[LPS2] C.--K. Li, Y.--T. Poon, and N.--S. Sze, Condition for the higher rank numerical range to be non--empty,
Linear and Multilinear Algebra 57, 365--368, 2009
\szs
[LS] C.--K. Li and N.--S. Sze, Canonical forms, higher--rank numerical ranges, totally isotropic subspaces,
and matrix equations, Proc. Amer. Math. Soc. 136, 3013--3023, 2008
\szs
[MM\.Z] K. Majgier, H. Maassen, and K. \.Zyczkowski, Protected subspaces in quantum information, Quantum Information
Processing 9, 343--367, 2010
\szs
[M] N. Mudalige, Higher Rank Numerical Ranges of Normal Operators, MSc thesis, U of
Guelph, 2010
\szs
[QD] J. F. Queir\'o and A. L. Duarte, Imbedding conditions for normal matrices, Linear Algebra Appl. 430, 1806--1811, 2009
\szs
[Wi] J. P. Williams, On compressions of matrices, J. London Math. Soc. (2) 3, 526--530, 1971
\szs
[Wo] H. Woerdeman, The higher--rank numerical range is convex, Linear and Multilinear Algebra 56, 65--67, 2008
\newpage\noindent
Author addresses:
\szs
John Holbrook\\
Dept of Mathematics and Statistics\\
University of Guelph\\
Guelph, Ontario, Canada N1G 2W1\\
jholbroo@uoguelph.ca
\szs
Nishan Mudalige\\
Dept of Mathematics and Statistics\\
York University\\
Toronto, Ontario, Canada\\
nishanm@yorku.ca
\szs
Rajesh Pereira\\
Dept of Mathematics and Statistics\\
University of Guelph\\
Guelph, Ontario, Canada N1G 2W1\\
pereirar@uoguelph.ca
\end{document}